\documentclass{amsart}
\newtheorem{Theo}{Theorem}
\begin{document}
\title{On Fabry's Gap Theorem}
\author[J.-C. Schlage-Puchta]{Jan-Christoph Schlage-Puchta}
\begin{abstract}
By combining Tur\'an's proof of Fabry's gap theorem with a gap theorem of P.
Sz\"usz we obtain a gap theorem which is more general then both these theorems.
\end{abstract}
\maketitle
Hadamard's classical gap theorem states that if $f(z)=\sum_{l=1}^\infty a_l
z^{k_l}$ is a power series with radius of convergence 1, and
$\frac{k_{l+1}}{k_l}\geq\theta>1$, the circle $|z|=1$ is the natural boundary
of $f$. The condition on the growth of the $k_l$ was lessened
by Fabry to
$k_l/l\rightarrow\infty$ without changing the conclusion(see \cite{Lan}).
Using his powersum
method, P. Tur\'an gave a simple proof of this theorem(see \cite{Tur} or
\cite{Tur2}, section 20). By a completely
elementary argument, P. Sz\"usz\cite{Sz} proved the following theorem:
\begin{Theo}
Suppose that $f(z)=\sum_{k=1}^\infty a_k z^k$ is a power series with radius of
convergence 1. Assume that there is a subsequence $k_l$ of indices, such that
$\limsup \sqrt[k_l]{a_{k_l}}=1$, and a constant $\theta>0$, such that $a_k=0$
for all $k\neq k_l$, $|k-k_l|<\theta k_l$. Then the circle $|z|=1$ is the
natural boundary of $f$.
\end{Theo}
Note that this theorem is stronger then Hadamard's gap theorem, but neither
implies nor is implied by Fabry's gap theorem. In this note we will combine
P. Tur\'an's proof of Fabry's gap theorem with the idea of P. Sz\"usz to obtain
the following theorem, which implies both gap theorems.
\begin{Theo}
Suppose that $f(z)=\sum_{k=1}^\infty a_k z^k$ is a power series with radius of
convergence 1. Assume that there is a subsequence $k_l$ of indices, such that
$\limsup \sqrt[k_l]{a_{k_l}}=1$, and a constant $\theta>0$, such that the
following holds true: If $N_l$ denotes the number of indices $k$ within the
interval $[(1-\theta)k_l, (1+\theta)k_l]$, such that $a_k\neq 0$, then
$N_l/k_l\rightarrow 0$. Then the circle $|z|=1$ is the natural boundary of $f$.
\end{Theo}
To obtain Fabry's gap theorem from this theorem, we can choose the sequence of
all nonvanishing coefficents, to deduce Sz\"usz' gap theorem, note that under
the assumptions of theorem 1 we have $N_l=1$ for any $l$.

To prove theorem 2, it suffices to show that $f$ is singular in 1. Assume that
$f$ was regular in 1. Then there was some $c>0$, such that $f$ is regular in
$\{|z-1|<c\}\cup\{|z|<1\}$. Choose $\delta<\theta/2$, and $\varphi>0$ so small,
that $f$ is regular within the domain $|z|<1+c/2, |\arg z| < 2\varphi$. Then we
will consider $f$ on the arc $\mathcal{C}=\{|z|=\delta, |\arg z|<\varphi\}$. First
we give an upper bound for $|f^{(k)}(z)|$ on $\mathcal{C}$. We may assume that
$\delta$ is so small, that around any point of $\mathcal{C}$, there is a disc with
radius $1+\delta$, such that $f$ is holomorphic within this disc. By the
standard integral estimate we obtain $|f^{(k)}(z)|<k!(1+\delta/2)^{-k}<k!/2$
for $k$ sufficently large.

To obtain a lower bound, we want to use Tur\'an's powersum technic. Let $k$ be
an index occuring in the sequence of our theorem and choose $m$ to be the
greatest integer such that $\left[\frac{m}{1-\delta}\right]\leq k$. Then we
have
\[
\frac{f^{(m)}(z)}{m!} = \sum_{\nu=0}^\infty {\nu\choose m} z^{\nu-m} a_\nu =
\sum_{\nu=m}^{(1+\delta)k} {\nu\choose m} z^{\nu-m} a_\nu +
\sum_{\nu>(1+\delta)k} {\nu\choose m} z^{\nu-m} a_\nu
\]
Now we restrict $z$ to be an element of $\mathcal{C}$.
We bound the second sum first. We have by the restriction $\nu>(1+\delta)k$
\[
\frac{{{\nu+1}\choose m}}{{\nu\choose m}}\delta < 2/3
\]
Together with the bound $a_\nu<(6/5)^\nu$, valid for $\nu$ sufficiently large, we see
that the second sum is $<1/2$.

Combining these estimates we obtain the following inequality:
\begin{equation}
\max_{z\in\mathcal{C}}\left|\sum_{\nu=m}^{(1+\delta)k} {\nu\choose m} z^{\nu-m} a_\nu
\right| < 1 
\end{equation}
The left hand side can be bounded from below using the continuos version of the
first main theorem in the powersum theory (see \cite{Tur2}, section 6):
\begin{Theo}
For $1\leq\nu\leq n$ let $\alpha_\nu$ be complex numbers, $\lambda_\nu$ real
numbers. Then we have for every $\mu>0$ the inequality
\[
\max\limits_{0\leq x\leq \mu}\left|\sum_{\nu=1}^n \alpha_\nu e^{i\lambda_\nu x}
\right| \geq \left(\frac{\mu}{4e\pi}\right)^n\max\limits_{0\leq x\leq 2\pi}
\left|\sum_{\nu=1}^n \alpha_\nu e^{i\lambda_\nu x}\right|
\]
\end{Theo}
Applying this to the left hand side of (1), we obtain the bound
\[
\max_{|z|=\delta}\left|\sum_{\nu=m}^{(1+\delta)k} {\nu\choose m} z^{\nu-m} a_\nu
\right| < \left(\frac{2e\pi}{\varphi}\right)^{N_l}
\]
The maximum of this sum is at least its quadratic mean, and the latter can be
computated using Parseval's equality, thus we get the bound
\[
\sum_{\nu=m}^{(1+\delta)k} {\nu\choose m}^2 \delta^{2(\nu-m)} |a_\nu|^2 < 
\left(\frac{2e\pi}{\varphi}\right)^{2N_l}
\]
Neglecting all terms on the left hand side except $\nu=k$, we get
\[
{k\choose m}^2 \delta^{2(k-m)} |a_k|^2 < \left(\frac{2e\pi}{\varphi}\right)^
{2N_l}
\]
By assumption we have $N_l=o(k_l)$, thus the right hand side is $e^{o(m)}$. The
same is true for $|a_k|$, and computing the binomial coefficient using
Stirling's
formula, we get the final contradiction
\[
\frac{1}{(1-\delta)^{2m}} = e^{o(m)}
\]
which finishes our proof.

\end{document}